\documentstyle[12pt]{amsart}

\theoremstyle{plain}
 \newtheorem{thm}{Theorem}[section]
 \newtheorem{prop}[thm]{Proposition}
 \newtheorem{lemma}[thm]{Lemma}
 
\theoremstyle{definition}
 \newtheorem{defn}{Definition}[section]

\theoremstyle{remark}

 \newtheorem{ack}{Acknowledgment}
      
\begin{document}
\title{Partial Gauss decomposition, \bf $U_q(\widehat{{{\frak{gl}}}(n-1)})\in
U_q(\widehat{{{\frak{gl}}}(n)}) $ and Zamolodchikov algebra}
\author{Jintai Ding}
\address{Jintai Ding,  RIMS, Kyoto University}
\maketitle

\begin{abstract}
We use the idea of partial Gauss decomposition to study structures related to 
$U_q(\widehat{{{\frak{gl}}}(n-1)})$ inside $U_q(\widehat{{{\frak{gl}}}(n)}) $.
This gives a description of $U_q(\widehat{{{\frak{gl}}}(n)})$ as an 
extension of $U_q(\widehat{{{\frak{gl}}}(n-1)})$ with  Zamolodchikov algebras, 
We explain the connection of this new realization with form factors. 
\end{abstract}

\section{ Introduction}

The affine
Kac-Moody algebra $\hat {\frak{g}}$ associated to a simple Lie algebra
$\frak{g}$
admits a natural realization as a central extension of the
corresponding loop algebra ${\frak{ g}} \otimes \Bbb C [t, t^{-1}] $.
Drinfeld gives a similar realization for $U_q(\frak{g})$, which 
is called Drinfeld realization\cite{D1}.   Faddeev, Reshetikhin and
Takhtajan  \cite{FRT}  Reshetikhin and Semenov-Tian-Shansk
\cite{RS} present  a  realization of
$U_q(\frak{g})$ to the quantum loop algebra $U_q(
{\frak{ g}} \otimes \Bbb [t, t^{-1}])$ using a solution
 of the Yang-Baxter equation depending on a parameter $z\in\Bbb C$
$$
R_{12}(z)R_{13}(zw) R_{23}(w)=R_{23}(w)R_{13}(zw) R_{12}(z),
$$
where $R(z)$ is a rational function of $z$ with values in
$\text{End}(\Bbb C^n\otimes \Bbb C^n)$.
An explicit identification 
between the two realizations of the quantum affine algebra $U_q(\hat
{\frak{g}})$ for the case ${\frak g}={\frak gl(n)}$ is established 
\cite{DF} by applying Gauss decomposition to the L-operators for the 
FRTS realization.

In this paper, we will use the idea of partial Gauss decomposition to 
study the structures related to 
$U_q(\widehat{{{\frak{gl}}}(n-1)})$ inside $U_q(\widehat{{{\frak{gl}}}(n)}) $.
We show  that $U_q(\widehat{{{\frak{gl}}}(n)})$ can be described as an 
extension of $U_q(\widehat{{{\frak{gl}}}(n-1)})$ with  Zamolodchikov algebras, 
where the   Zamolodchikov algebras can be interpreted as intertwiner for 
 $U_q(\widehat{{{\frak{gl}}}(n-1)})$.  The Zamolodchikov algebras are
use to derive structures related to form factors and related structures, 

This paper is to present the method of the partial Gauss decomposition 
to find new structures hidden inside the affine quantum groups, which  is 
related to many aspects of the theory of affine quantum groups
\cite{Di2}\cite{MJ}\cite{M} and physics\cite{Sm}. 
The meaning of this method can  be explained using the method 
of the twisting of Drinfeld \cite{D2}\cite{KT}.

\section{ Quantum algebra $U_q(\hat{\frak{gl}}(n-1))
\in U_q(\hat{\frak{gl}}(n))$ and partial Gauss decomposition}

Let $V$ be $\Bbb C^n$ with a fixed basis $e_i, i=1,..,n$
and $E_ij$ be the standard basis of End$(\Bbb C^n)$ dependent 
on $e_i$. 
Let  $R(z)$ be an element of End$(\Bbb C^n\otimes \Bbb C^n)$ defined by 
\begin{align*}
R(z) &= \sum^n_{i=1}E_{ii}\otimes E_{ii} + \sum^n\Sb i\neq
j\\ i,j=1\endSb E_{ii}\otimes E_{jj}\frac{z-1}{q^{-1}z
-q}\cr 
&+ \sum^n\Sb i>j\\ i,j=1\endSb E_{ij}\otimes
E_{ij}\frac{(q^{-1}-q)}{zq^{-1}-q} + \sum^n\Sb i<j\\
i,j=1\endSb E_{ij}\otimes E_{ji} \frac{z(q^{-1}-q)}{zq^{-1}-q}
\end{align*}
where $q,z$ are formal variables. Then 
$R(z)$ satisfies the  Yang-Baxter equation 
 and $R$ is unitary, namely 
$$
R_{21}(z)^{-1} = R({z^{-1}}),
$$
where $R_{21}(z)= PR_{12}(z)P$, where $P$ is the operator
permuting the two components $V\otimes V$. Here $V=\Bbb C^n$.

Faddeev, Reshetikhin and Takhtajan defined a Hopf  algebra using 
 $R(z)$, which satisfies Yang-Baxter equation.  Reshetikhin and
Semenov-Tian-Shansky obtained a central extension of this algebra.
The algebra defined with the $R(z)$ above 
is isomorphic $U_q(\hat{\frak{gl}}(n))$.  
The central extension is incorporated in shifts of the parameter $z$
in $R(z)$.

\begin{defn} 
$U_q(\widehat{\frak{gl}}(n))$ is an associative algebra
with generators $\{l^{\pm}_{ij}[\mp m], m\in{\bold Z_+ \setminus  0},
l^+_{ij}[0], l^-_{ji}[0], 1\leq j\leq i\leq n\}$.  Let
$l^{\pm}_{ij}(z) = \sum\limits^{\infty}_{m=0} l^{\pm}_{ij}[\pm
m] z^{\pm m}$, where $l^+_{ij}[0] = l^-_{ji}[0] = 0$, for $1\leq
j>i\leq n$.  Let $L^{\pm}(z) = (l^{\pm}_{ij}(z))^n_{i,j=1}$.
Then  the defining  relations are the following:
 $$
l^+_{ii}[0] l^-_{ii}[0] = l_{ii}[0] l^+_{ii}[0]=1,
$$
\begin{align*}
R(\frac{z}{w}) L^{\pm}_1(z) L^{\pm}_2(w) &= L^{\pm}_2(w) L^{\pm}_1(z)
R(\frac{z}{w}),\\
R(\frac{z_-}{w_+}) L^+_1(z) L^-_2(w) &= L^-_2(w) L^+_1(z)
R(\frac{z_+}{w_-}),
\end{align*}
where $z_{\pm} = zq^{\pm  \frac{c}{2}} $.
The expansion direction of
$R(\frac{z}{w})$ are chosen to be in $\frac{z}{w}$ or $\frac{w}{z}$
respectively. \cite{DF}
\end{defn}

The  Hopf algebra is given by:
\begin{align*}
\Delta L^{\pm}(z) = L^{\pm}(zq^{\pm(1\otimes\frac{c}{2})})
&\dot{\otimes} L^{\pm}(zq^{\mp(\frac{c}{2}\otimes 1)})\\ 
\text{or}\qquad \Delta(l^{\pm}_{ij}(z)) = \sum^n_{k=1}
l^{\pm}_{ik}(zq^{\pm(1\otimes\frac{c}{2})}) &\otimes
l^\pm_{kj}(zq^{\mp(\frac{c}{2}\otimes 1)}),
\end{align*}

and its antipode is
$$
S(L^{\pm}(z)) = L^{\pm}(z)^{-1}.
$$
The invertibilty of $L^{\pm}(z)$ follows from the properties  that
  $l^{\pm}_{ii}$ are invertible and $L^{\pm}(0)$ are upper triangular and
lower triangular, respectively. 

$L^{\pm}(z)$ have the following unique decompositions:
\begin{align*}
L^{\pm}(z) = &\pmatrix 1& & &  0\\
 e^{\pm}_{2,1}(z) &\ddots & &\\ e^{\pm}_{3,1}(z)
&\ddots &\ddots\\ \vdots &\ddots &\ddots &\ddots\\ e^{\pm}_{n,1}(z)
&\hdots &e^{\pm}_{n,n-1}(z) &e^{\pm}_{n-1,n}(z) &1\endpmatrix
\pmatrix k^{\pm}_1(z) & & &  0 \\ &\ddots\\  & &\ddots\\ &\\0 & &
&k^{\pm}_n(z)\endpmatrix\times \tag 3.22 \\  \\
&\pmatrix 1 & f^{\pm}_{1,2}(z) & f^{\pm}_{1,3}(z) &\hdots
&f^{\pm}_{1,n}(z)\\ &\\ &\ddots &\ddots &\ddots &\vdots\\ & & &
&f^{\pm}_{n-1,n}(z)\\ &\\0  & & & &1\endpmatrix,
\end{align*}
which is used to establish the isomorphism between Drinfeld realizations 
of $U_q(\widehat{{{\frak{gl}}}(n)})$ and its FRTS realization. 

Similarly we have the following partial Gauss decomposition:
\begin{prop} 
The operator $L^\pm(z)$ can be uniquely decomposed as 
\begin{align*}
L^{\pm}(z) = &\pmatrix I& 0\\
 e^{\pm}(z) &1 \endpmatrix
\pmatrix K^{\pm}(z) & 0 \\ 0 &k^{\pm}(z)\endpmatrix
\pmatrix I & f^{\pm}(z)\\ 0 &1\endpmatrix,
\end{align*}
where $K^{\pm}(z)$ and $k^{\pm}(z)$ are $n-1\times n-1$  
invertible matrix operators,
 $e^{\pm}(z)$ is a size $n-1$ column and $f^{\pm}(z)$ is a size $n-1$ arrow.  
\end{prop}
Because $K^{\pm}(z)$ are invertible,
the elements $e^{\pm}(z)$, $f^{\pm}(z)$ and
$k^{\pm}(z)$  are uniquely expressed in terms of the matrix coefficients of 
$L^\pm(z)$.  

We have: 
\begin{prop}
The algebra generated by entries of  operator matrixes  $K^{\pm}(z)$
is  $U_q(\hat {\frak gl}(n-1))$. 
\end{prop}

We will follow th step as in the case of \cite{DF} to find out the 
complete commutation relations for the operators in the above decomposition. 
For the calculation, we need the following formulas: 
\begin{align*}
L^{\pm}(z) &= \pmatrix K^{\pm}(z) & K^{\pm}(z)f^{\pm}(z)\\
\\e^{\pm}(z)K^{\pm}(z) & (k^{\pm}(z) + e^{\pm}(z)
K^{\pm}(z)f^{\pm}(z)) 
\endpmatrix, \\   \\
&L_1(z)L_2(w) = \\
&\pmatrix K^{}(z)K^{}(w) &K^{}(z) K^{}(w)f^{}(w)
&K^{}(z)f^{}(z)K^{}(w) &K^{}(z)f^{}(z) K^{}(w)f^{}(w)
\\  \\K^{}(z) e^{}(w)K^{}(w) &K^{}(z) D^{}(z)
& K^{}(z)f^{}(z)
 e^{}(w)K^{}(w) &K^{}(z)f^{}(z) D^{}(w) \\ \\
e^{}(z)K^{}(z)K^{}(w) &e^{}(z)K^{}(z)K^{}(w)f^{}(w)
&D^{}(z)K^{}(w) &D^{}(z)
(K^{}(w)f^{}(w))
\\  \\
e^{}(z)K^{}(z) e^{}(w)K^{}(w) &e^{}(z)K^{}(z)D^{}(w)
& D^{}(z)  e^{}(w)K^{}(w) & D^{}(z)(D^{}(w))
K^{}(w)f^{}(w) 
\endpmatrix, \\   \\
&(L^{\pm}(z))^{-1} = \pmatrix K^{\pm}(z)^{-1} +
f^{\pm}(z)k^{\pm}(z)^{-1}e^{\pm}(z) &
-f^{\pm}(z)k^{\pm}(z)^{-1}\\ \\ -k^{\pm}(z)^{-1}e^{\pm}(z)
&k^{\pm}(z)^{-1}\endpmatrix,  \\ \\
&R_{12}(\frac{z}{w}) = \pmatrix \bar R(z/w) & 0 & 0 & 0\\ 0 &
\frac{z-w}{zq^{-1}-wq}A &\frac{-z(q-q^{-1})}{zq^{-1}-wq}B &0\\ 0 &
\frac{-w(q-q^{-1})}{zq^{-1}-wq}C &\frac{z-w}{zq^{-1}-wq}D &0\\ 0 & 0 & 0
& 1\endpmatrix \\ \\ 
&R_{12}(\frac{z}{w})^{-1}= 
R_{21}(\frac{w}{z}) = \pmatrix \bar R_{21}(w/z) & 0 & 0 & 0\\ 0 &
\frac{z-w}{zq-wq^{-1}}A &\frac{-z(q-q^{-1})}{wq^{-1}-zq}B &0\\ 0 &
\frac{-w(q-q^{-1})}{wq^{-1}-zq}C &\frac{z-w}{zq-wq^{-1}}D &0\\ 0 & 0 & 0
& 1\endpmatrix \\ \\
\end{align*}
\begin{align*}
&(L_1(w))^{-1} =\\
&\pmatrix \ast & 0 & -f(w)k(w)^{-1} & 0\\
0 & \ast & 0 & -f(w)k(w)^{-1}\\
-k(w)^{-1}e(w) & 0 & k(w)^{-1} & 0\\
0 & -k(w)^{-1}e(w) & 0 & k(w)^{-1}\endpmatrix, 
\end{align*}

where 
$A= \sum_{i\neq n} E_{ii}\otimes E_{nn}$, $D= 
\sum \Sb {i\neq n} \endSb E_{nn}\otimes E_{ii}$
$C=\sum\Sb j<n+1 \endSb E_{nj}\otimes
E_{jn}$, $B= 
\sum \Sb j<n+1 \endSb
E_{jn}\otimes E_{nj}$, 
$\bar R(z)$ is the R-matrix restricted to the subapce 
$V'\otimes V'$, $V'$ is generated on the subspace generated by 
$e_i, i=1,..,n-1$, 
$D^{\pm}(z)= (k^{\pm}(z) + e^{\pm}(z)
K^{\pm}(z)f^{\pm}(z))$ ;  and

$$L^{\pm}_1(w)^{-1}R_{21}(\frac{z}{w}) L^{\pm}_2(z) =
L^{\pm}_2(z) R_{21}(\frac{z}{w})L^{\pm}_1(w)^{-1}$$ 
$$  L^-_1(w)^{-1} R_{21}(\frac{z^{+}}{w^{-}})L^+_2(z) =
L^+_2(z) R_{21}(\frac{z_-}{w_+})L^-_1(w)^{-1}$$ 
$$R_{21}(\frac{z_-}{w_+})L^-_2(z) L^+_1(w) =
L^+_1(w)L^-_2(z) R_{21}(\frac{z_+}{w_-})$$ 
$$L^+_1(w)^{-1} R_{21}(\frac{z_-}{w_+})L^-_2(z) = L^-_2(z)
R_{21}(\frac{z_+}{w_-})L^+_1(w) $$
$$L^{\pm}_2(z)^{-1}(L^{\pm}_1(w))^{-1} R(\frac{z}{w}) =
R_{21}(\frac{z}{w}) (L^{\pm}_1(w))^{-1}(L^{\pm}_2(z))^{-1}$$
$$L^+_2(z)^{-1} L^-_1(w)^{-1}R_{21}(\frac{z_+}{w_-}) =
R_{21}(\frac{z_-}{w_+}) (L^-_1(w))^{-1}(L^+_2(z))^{-1}$$

Using the same calculation technique as in \cite{DF}, we  have:
\begin{lemma}

\begin{align*}
\bar R(z/w)K_1^{\pm}(z)K_2^{\pm}(w) &= K_2^{\pm}(w) K_1^{\pm}(z)\bar R(z/w)\\
k^{\pm}(z) k^{\pm}(w) &= k^{\pm}(w) k^{\pm}(z)\\
\bar R(z_+/w_-)K_1^+(z)K_2^-(w) &= K_2^-(w)K_1^+(z)\bar R(z_-/w_+)\\
k^+(z)k^-(w) &= k^-(w)k^+(w)\\
k^{\pm}(z)k^{\pm}(w) &= k^{\pm}(w)k^{\pm}(z)\\
\frac {z_{\mp}q^{-1}-w_{\pm}q}{z_{\mp }-w_{\pm}}
k^{\mp}(w)^{-1}K^{\pm}(z) &=
K^{\pm}(z)k^{\mp}(w)\frac {z_{\pm}q^{-1}-w_{\mp}q}{z_{\pm }-w_{\mp}},
\end{align*}

$$
K_1^\pm(z) E_2(w) = \frac{zq^{\mp\frac{c}{2}-1} - wq}{zq^{\mp\frac{c}{2}} -
w} E_2(w)\bar R(zq^{\mp\frac{c}{2}}/w)
 K_1^\pm(z),$$
$$
K_1^pm(z)\bar R(zq^{\pm\frac{c}{2}}/w) F_2(w)
 = \frac{zq^{\pm\frac{c}{2}} - w}{zq^{\pm\frac{c}{2}-1} -
wq} F_2(w) K_1^\pm(z),$$
$$k^\pm(z) E(w) = \frac{zq^{\mp\frac{c}{2}+1} -
wq^{-1}}{zq^{\mp\frac{c}{2}} - w} E(w) k^\pm(z),$$
$$k^\pm(z) F(w) = \frac{zq^{\pm\frac{c}{2}} - w}{zq^{\pm\frac{c}{2}+1} -
wq^{-1}} F(w) k^\pm(z),$$

$$({z-wq^2})E_1(z) E_2(w)R(z/w) = ({zq^2-w}) E_2(w) E_1(z),$$
$$({zq^2-w})F(z) F(w) =R(z/w) ({z-wq^2}) F(w) F(z),$$
\begin{multline*}
 E_2(z)(F_1(w))-
F_1(w)E_2(z)  = \\ (q-q^{-1}) \left( \delta \left(\frac{w}{z}
q^c \right) k^-(w q^{\frac{c}{2}}) K^-(w q^{\frac{c}{2}})^{-1} - \delta
\left(\frac{w}{z} q^{-c} \right) k^+(w q^{-\frac{c}{2}}) K^+(w
q^{-\frac{c}{2}})^{-1} \right),
\end{multline*}
where  $E(z)=e^+(zq^{\frac{c}{2}})-e^-(zq^{\frac{-c}{2}})$, 
       $F(z)=f^+(zq^{\frac{-c}{2}})-f^-(zq^{\frac{c}{2}})$  and  
       $$\delta(x) = \sum_{m\in\Bbb Z} x^m.$$
\end{lemma}

 The algebra generated by E(z) or F(z) gives a realization of the 
Zamolodchikov algebra, the formulation above is basically the same as in 
\cite{Di1}, where we study the Hopf algebra extension of 
Zamolodchikov algebras. On the other hand, we can reformulate the 
definition of $U_q(\hat {\frak gl(n)})$ using the relations above. 

\begin{defn}
Let ZUR(n) be an algebra generated by matrix operators $ E(z),$ $F(z)$ 
$K^{\pm}(z)$ and $k^\pm(z)$ associated with  the vector space $V=\Bbb C^{n-1}$ 
respectively to $V^*,V, V\otimes V^*$ and a one dimesional space $\Bbb C$. 
The commutation relations are defined as in the lemma above. 
\end{defn}

Then 
\begin{thm}
 ZUR(n) is isomorphic to $U_q(\hat {\frak gl}(n))$. 
\end{thm}

The proof follows form the above lemma and the similar argument in 
\cite{DF}. 

From the point of view of \cite{Di1}, we can similarly to give a new 
Hopf algebra structure to this formulation using the 
similar formulas. 
The important point is that from the definition we can see that 
$E(z)$ and $F(z)$ is nothing but intertwiner for the 
affine algebra $U_q(\hat {\frak gl}(n-1))$ generated by the 
operators $K^\pm(z)(k^\pm(z))^{-1}$.   The last formula of 
the commutation relations implies the constructions like in \cite{M}
\cite{Sm}.

Let $\bar E(z)= E(z) K^-(zq^{c/2}k^-(zq^{c/2})$, then 
\begin{prop}
\begin{align*}
E_2(z)\frac {(1-z/wq^{-c})R(z/w)(z-wq^2)}{zq^2-w}(F_1(w))-
F_1(w)E_2(z)  = \\ 
(q-q^{-1}) (1-q^{-2c}) \delta(\frac w z q^{-c})
k^+(w q^{-\frac{c}{2}}) K^+(w
q^{-\frac{c}{2}})^{-1}k^-(w q^{\frac{c}{2}})^{-1} K^-(w q^{\frac{c}{2}}),
\end{align*}
$$({zq^2}-w)E_1(z) E_2(w) = ({z-wq^2}) E_2(w) E_1(z)R(z/w),$$
$$({zq^2-w})F(z)_1 F_2(w) =R(z/w) ({z-wq^2}) F_2(w) F_1(z),$$
$$
(1-z/wq^{c} E_2(z)R(z/w)\frac {z-wq^2}{zq^2-w}(F_1(w))-
F_1(w)E_2(z)  = (q-q^{-1}) (1-q^{2c}) \delta(\frac w z q^{c}). 
$$
\end{prop}

The first one of the formulas above coincides with the spinor constructions of 
affine quantum groups in \cite{Di2}. 

This last three formula above 
says that these operators generate an algebra almost the same as  
the Zamolodchikov-Faddeev algebra used to describe 
the theory of form factors\cite{Sm}. Similarly 
we can also define a new operator 
$\bar F(z)=( k^+(w q^{-\frac{c}{2}})^{-1} K^+(w
q^{-\frac{c}{2}})F(z)$. This operator with $E(z)$ generate another  algebra 
similar to the definition above. 
From  the point view  of intertwiners as in \cite{MJ}, those operator
can give a complete theory of form factors, where one copy of the 
algebra is explained as the the Zamolodchikov-Faddeev algebra to define the 
model and the other one is explained as local operators, which commutes
with the first algebra up to certain functions. In a subsequent paper. 
we will apply the same method to Yangian, and the elliptic algebra \cite{LKP}
\cite{F} We will give the complete details to the descriptions of 
more general Zamolodchikov-Faddeev type of algebras, whose degeneration
gives us the corresponding results in this paper.

\begin{ack}
We would like to thank B. Feigin and T. Miwa for useful discussions. 
\end{ack}

\end{document}